\documentclass[10pt]{amsart}
\usepackage{geometry}                
\usepackage{graphicx}
\usepackage{amssymb}
\usepackage{epstopdf}
\newtheorem{theorem}{Theorem}
\newtheorem{corollary}[theorem]{Corollary}
\newtheorem{lemma}[theorem]{Lemma}

\def\eproof{\rule{3mm}{3mm}\newline}

\title{Hyperbolic volume of n-manifolds with geodesic boundary and orthospectra}
\author{Martin Bridgeman}
\address{Math Dept., Boston College, Chestnut Hill, Ma 02167}
\email{bridgem@bc.edu}
\author{Jeremy Kahn}
\address{Math Dept., Stony Brook University, Stony Brook, NY 11794}
\email{kahn@math.sunysb.edu}
\date{Dec 17, 2009}                                            
\thanks{Bridgeman was partially supported by NSF grant DMS-0707116
and Kahn was partially supported by NSF grants DMS-0905812}
\begin{document}
\maketitle
\begin{abstract}
In this paper we describe a function $F_n:{\bf R}_+ \rightarrow {\bf R}_{+}$ such that for any  hyperbolic n-manifold $M$ with totally geodesic boundary $\partial M \neq \emptyset$, the volume of $M$ is equal to the sum of the values of $F_n$  on the {\em orthospectrum} of $M$. We derive an integral formula for $F_n$ in terms of  elementary functions. We use this to give a lower bound for the volume of a hyperbolic n-manifold with totally geodesic boundary in terms of the area of the boundary.
\end{abstract}

 \section{Introduction}

We let $M$ be a compact hyperbolic $n-$manifold with non-empty totally geodesic boundary. An orthogeodesic for $M$ is a geodesic arc  with endpoints in $\partial M$ and perpendicular to $\partial M$ at the endpoints. Then the {\em orthospectrum}  $\Lambda_M$  of $M$ is the set (with multiplicities) of lengths of orthogeodesics. Let $DM$ be the closed manifold obtained by doubling $M$ along the boundary, then $DM$ is a closed hyperbolic manifold and therefore the set of closed geodesics of $DM$ is countable. As the orthogeodesics of $M$ correspond to a subset of the closed geodesics of $DM$, the set of   orthogeodesics of $M$ is also countable and therefore $\Lambda_M$ is also. By decomposing the unit tangent bundle of $M$ we obtain the following theorem.

\begin{theorem}
Given $n \geq 2$ there exists a continuous monotonically decreasing function $F_n:{\bf R}_+ \rightarrow {\bf R}_+$ such that if $M$ is a compact hyperbolic $n-$manifold with non-empty totally geodesic boundary, then
$$\mbox{Vol}(M) = \sum_{l \in \Lambda_M} F_n(l).$$
\end{theorem}

We give an integral formula for $F_n$  over the unit interval of an elementary function and show show that  $F_n$  satisfies $\lim_{l \rightarrow 0^+} l^{n-2}F_n(l)= K_n > 0$. Therefore for $n \geq 3$, if a hyperbolic n-manifold with geodesic boundary has a short orthogeodesic, it has large volume. Conversely, if it doesn't have a short orthogeodesic, the boundary has a large embedded neighborhood and therefore large volume. Using this we prove the following theorem.

\begin{theorem}
For $n \geq 3$, there exists a monotonically increasing function $H_n:{\bf R}_+ \rightarrow{\bf R}_+$ and a constant $C_n > 0$ such that if $M$ is a hyperbolic n-manifold with totally geodesic boundary of area $A$ then 
$$\mbox{Vol}_n(M) \geq H_n(A) \geq C_n.A^{\left(\frac{n-2}{n-1}\right)}$$
\label{vol}
\end{theorem}

\section{Decomposition via orthogeodesics}
In this section we define $F_n$ and show that it gives a decomposition of the volume of $M$ as described in Theorem 1.

We let $C_M$ be the set of orthogeodesics of $M$ and write $C_M = \{ \alpha_i \}_{i\in I}$. We further let  $l_i = \mbox{Length}(\alpha_i)$.

Let $v \in T_{1}(M)$ and let $\alpha_v$ be the maximal length geodesic arc in $M$ tangent to $v$. We define
$$T^f_1(M) = \{ v \in T_1(M)\left| \  \partial \alpha_v \in \partial M \}\right..$$

If  $v \in T_1^f(M)$ then $\alpha_v$ is a closed geodesic arc with endpoints in $\partial M$ intersecting $\partial M$ transversely. We define $\sim$ on $T^f_1(M)$ by letting $v \sim w$ if $\alpha_v$ is homotopic to $\alpha_w$ in $M$ rel boundary $\partial M$.
Let $v_i \in T_1(M)$ be such that $v_i$ is tangent to the orthogeodesic $\alpha_i$. Then obviously $\alpha_{v_i} = \alpha_i$ and we let  $D_i = [v_i]$. We will show that the $D_i$ are exactly the equivalence classes of $\sim$.

We consider the universal cover $\tilde{M}$ of $M$ in ${\bf H}^n$. Then $\partial{\tilde{M}}$ is a collection of disjoint planes $P_i$ bounding disjoint hyperbolic open half spaces $H_i$ such that $\tilde{M} = {\bf H}^3 - (\cup_i H_i)$.

Given a $v \in T_1^f(M)$ we lift $\alpha_v$ to a geodesic arc $\tilde{\alpha}_v$ in $\tilde{M}$. Then $\tilde{\alpha}_v$ has endpoints in two disjoint components $P_i, P_j$ of $\partial{\tilde{M}}$. As $P_i, P_j$ are disjoint, there is a unique perpendicular $\tilde{\beta}$ between them in $\tilde{M}$. Then $\tilde{\alpha}_v$ is homotopic to $\tilde{\beta}$ in $\tilde{M}$ rel boundary. We let $\beta$ be the geodesic arc obtained by projecting $\tilde{\beta}$ down to $M$. Then $\beta$ is a geodesic arc with endpoints in $\partial M$ and perpendicular to 
$\partial M$. Therefore $\beta$ is an orthogeodesic of $M$ and therefore $\beta = \alpha_k$ for some $k$. Also the homotopy between $\tilde{\alpha}_v$ and $\beta$ in $\tilde{M}$ rel boundary, descends to a homotopy between $\alpha_v$ and $\alpha_k$. Therefore if we take $v_k$ be a tangent vector to $\alpha_k$ then $[v] = [v_k] = D_k$.

Now to show that $D_i \cap D_j = \emptyset$ for $i \neq j$, we note that if $v_i \sim v_j$ then $\alpha_i$ is homotopic to $\alpha_j$ rel boundary.We lift this homotopy to obtain a homotopy between lifts $\tilde{\alpha}_i$ and $\tilde{\alpha}_j$ in $\tilde{M}$ rel boundary. Let $P_k, P_l$ be the components of  $\partial \tilde{M}$ joined by $\tilde{\alpha}_i$. Then $\tilde{\alpha}_i$ is the unique  perpendicular between $P_k$ and  $P_l$.
As $\tilde{\alpha}_j$ is homotopic rel boundary to $\tilde{\alpha}_i$, it must also  connect $P_k, P_l$ and  is also the unique  perpendicular between $P_k$ and  $P_l$. Thus $\tilde{\alpha}_j = \tilde{\alpha}_j$ adn therefore $\alpha_i = \alpha_j$.

By the ergodicity of geodesic flow on the double $DM$, almost every $\alpha_v$ must have both endpoints in $\partial M$. Therefore $T_1^f(M)$ is of full measure in $T_1(M)$. Therefore integrating over the fibers, we have

$$\Omega_M(T_1(M)) = \mbox{Vol}({\bf S}^{n-1}).\mbox{Vol}(M) = \sum_{i \in I}  \Omega_M(D_i)$$
giving
$$\mbox{Vol}(M) = \frac{1}{\mbox{Vol}({\bf S}^{n-1})}\sum_{i \in I}  \Omega_M(D_i)$$
We now show that $\Omega_M(D_i)$ depends only on $l_i$. 
Let $p:T_1(\tilde{M}) \rightarrow T_1(M)$ be the covering map associated to the covering $\pi:\tilde{M} \rightarrow M$. We let $\Omega$ be the standard volume measure on $T_1({\bf H}^n)$. Then $p$ is a local isometry between $\Omega$ and $\Omega_M$.

Given $v_i$ a tangent vector to $\alpha_i$, we lift to $\tilde{\alpha}_i$ in $\tilde{M} \subset {\bf H}^n$. Then $\tilde{\alpha}_i$ is the unique perpendicular between the planes $P_j, P_k$. We let $\tilde{D}_i$ be the set of tangent vectors in $T_1(\tilde{M})$ tangent to a geodesic arc with endpoints on $P_j, P_k$.

 If $v \in D_i$ then $\alpha_v$ is homotopic to $\alpha_i$ and therefore lifting the homotopy, it lifts to a geodesic arc $\tilde{\alpha}_v$ with endpoints in  $P_j, P_k$. Therefore we can lift any point of $D_i$ to a point of $\tilde{D}_i$ and the lift is a local isometry between $T_1(M)$ and $T_1(\tilde{M})$.  Also by projecting back down to $M$, every point of $\tilde{D}_i$ is a lift of a point of $D_i$. Therefore $p$ restricts to a covering map from $\tilde{D}_i$ to $D_i$. To show it is a homeomorphism, if $g$ is a covering transformation for the covering  $\pi:\tilde{M} \rightarrow M$ that sends $\tilde{v} \in \tilde{D}_i$ to $\tilde{w} \in \tilde{D}_i$, then $g$ must send the pair of boundary components $P_j, P_k$ to themselves (possibly switching). Therefore $g$ must preserve the perpendicular $\tilde{\alpha}_i$ and fix (at least) the center point of of $\tilde{\alpha}_i$. As
covering transformations have no fixed points this is a contradiction. Therefore $\tilde{D}_i$ is an isometric lift of $D_i$. Therefore $\Omega_{M}(D_i) = \Omega(\tilde{D}_i)$.

We now take the upper half space model for ${\bf H}^n$ and denote the planes $P_i$ by the disks in ${\bf R}^{n-1}$ bounded by the half space $H_i$. We define 
$$D = \{ x \in {\bf R}^{n-1} \left|  \ |x| \leq 1\}\right. \mbox{ and }  D_{a} = \{ x \in {\bf R}^{n-1} \left|  \ |x| \geq a \}\right.$$
and the associated half-spaces by  $H, H_a$ respectively.  For each $l_i$, we define $a_i = e^{l_i}$. Then by an isometry  we can  map the half-spaces $H_j, H_k$ to half-spaces $H, H_{a_i}$.
We let $Q_a$ be the set of tangent vectors in $T_1({\bf H}^n-(H\cup H_a))$ tangent to a geodesic arc with endpoints in $D, D_a$. Then we have  $\Omega(Q_{a_i}) = \Omega(\tilde{D}_i)$ and we define  $$F_n(l) = \frac{\Omega(Q_{a_i})}{\mbox{Vol}({\bf S}^{n-1})}.$$
Then
$$\mbox{Vol}(M) = \frac{1}{\mbox{Vol}({\bf S}^{n-1})}\sum_{i \in I}  \Omega_M(D_i) = \frac{1}{\mbox{Vol}({\bf S}^{n-1})}\sum_{i \in I}  \Omega(\tilde{D}_i) = \sum_{i \in I} F_n(l_i)$$

\section{Integral formula}
We consider the upper half space model of ${\bf H}^n$. As described above we have disks
$$D = \{ x \in {\bf R}^{n-1} \left|  \ |x| \leq 1\}\right. \mbox{ and }  D_{a} = \{ x \in {\bf R}^{n-1} \left|  \ |x| \geq a\}\right.$$ bounding planes $P, P_a$ respectively. Then $Q_a$ is the set of tangent vectors tangent to a geodesic arc with endpoints on $P_0, P_a$. 

We let $\Omega$ be the volume form on $T_1({\bf H}^n)$. Let $ v \in T_1({\bf H}^n)$ and define $g_v$ to be the directed geodesic with tangent vector $v$. Then $g_v$ can be parameterized by the ordered pair of endpoints $(x,y) \in {\bf R}^{n-1}\times  {\bf R}^{n-1}$. Furthermore if $v$ has basepoint $p$ then $p$ is a unique signed hyperbolic distance from the highest point (in the upper half-space model) of $g_v$. We can therefore uniquely parameterize $T_1({\bf H}^n)$ by triples $(x,y,t) \in {\bf R}^{n-1}\times  {\bf R}^{n-1}\times {\bf R}$ where $(x,y)$ is the directed endpoints of $g_v$ and $t$ is the signed distance from the highest point of $g_v$.

In terms of this parametrization we have
$$d\Omega = \frac{2dV(x)dV(y)dt}{|x-y|^{2n-2}}$$
where $dV(x) = dx_1dx_2\ldots dx_{n-1}$, for $x = (x_1, x_2,\ldots,x_{n-1}) \in {\bf R}^{n-1}$.

Now if $v \in Q_{a}$ then $g_v$ has endpoints in $D \times D_a$ or $D_a\times D$. We define the map 
$L_a:(D \times D_a)\cup(D_a \times D) \rightarrow {\bf R}_+$ given by letting $L_a(x,y)$ equal the length of the segment of the geodesic with endpoints $x,y$ between $P_0$ and $P_a$.

Then integrating along the $t$ direction we have
$$F_n(l) = \frac{1}{\mbox{Vol}({\bf S}^{n-1})} \int_{Q_a} d\Omega =  \frac{1}{\mbox{Vol}({\bf S}^{n-1})} \int_{Q_a} \frac{2dV(x)dV(y)dt}{|x-y|^{2n-2}}$$
Giving 
$$F_n(l) = \frac{2}{\mbox{Vol}({\bf S}^{n-1})} \int_{(D \times D_a)\cup(D_a \times D)}
\frac{L_a(x,y)dV(x)dV(y)}{|x-y|^{2n-2}}$$

As $L_a(x,y) = L_a(y,x)$ we have
$$F_n(l) = \frac{4}{\mbox{Vol}({\bf S}^{n-1})} \int_{D \times D_a}\frac{L_a(x,y)dV(x)dV(y)}{|x-y|^{2n-2}} = \frac{4}{\mbox{Vol}({\bf S}^{n-1})} \int_{|x| \leq 1}\int_{|y| \geq a} \frac{L_a(x,y)dV(x)dV(y)}{|x-y|^{2n-2}}$$

We first note the following fact from prior work.
\begin{lemma}{(Bridgeman \cite{B09})}
If $x,y \in {\bf R}$ then 
$$L_a(x,y) = \frac{1}{2}\log\left(\frac{(y^2-1)(x^2-a^2)}{(y^2-a^2)(x^2-1)}\right)$$
\end{lemma}
It follows that
$$F_2(l) = \frac{1}{\pi} \int_{-1}^1\int_{a}^{\infty} \frac{\log\left(\frac{(y^2-1)(x^2-a^2)}{(y^2-a^2)(x^2-1)}\right)dxdy}{(x-y)^{2}}.$$

The above integral can be explicitly calculated yielding the following;
\begin{theorem}{(Bridgeman \cite{B09})}
For $n = 2$  
$$F_{2}(l) = \frac{4}{\pi}{\mathcal L}\left(\frac{1}{\cosh^2\frac{l}{2}}\right)$$
where ${\mathcal L}$ is the Rogers $L-$function.
\end{theorem}

\section{Higher dimensional case}
We now consider the case of $n \geq 3$. Let $(x,y) \in D\times D_a$.   
We let $z = y-x$ with $u = z/|z|$. Also we define $w$ to be the point on the line $l$ joining $x, y$ in ${\bf R}^{n-1}$ nearest to the origin in the Euclidean metric on ${\bf R}^{n-1}$. Then $w = r.v$  where $v \in {\bf S}^{n-2}$ and $r > 0$.  Then $u, v \in  {\bf S}^{n-2}$ are perpendicular and
$$x =  s.u+r.v  \qquad y =  t.u + r.v$$
for unique $s,t \in {\bf R}$.

\begin{figure}[htbp] 
   \centering
   \includegraphics[width=4in]{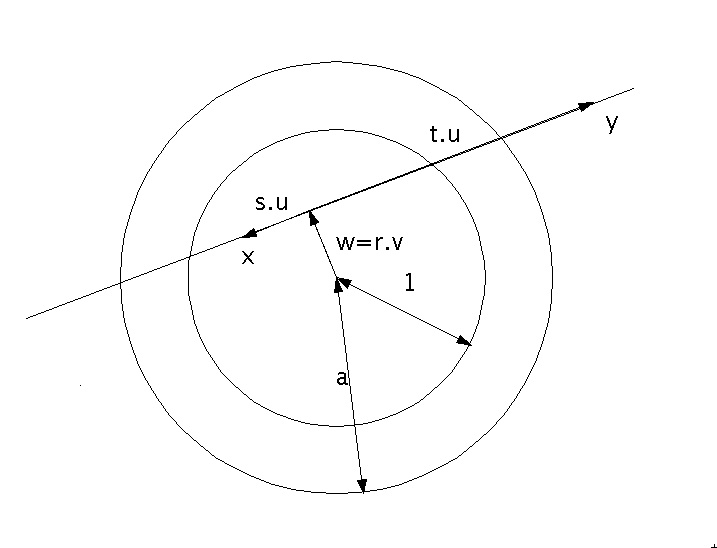} 
   \caption{Parametrization by $r,s,t,u,v$}
\end{figure}
We reparametrise the integral in terms of $r,s,t, u,v$. We note that $(u,v)$ parametrize the unit tangent bundle  $T_1({\bf S}^{n-2})$. 

To calculate the Jacobian, we first change variables from $(x,y)$ to $(x, z)$ where $z = y-x$. Then the Jacobian is trivially equal to one. As $z = y-x = (t-s)u$, then from standard polar coordinates we have that $$dV(z) = (t-s)^{n-2}d\Omega(u)d(t-s).$$
 where $d\Omega(u)$ is the volume form on the unit sphere  ${\bf S}^{n-2}$ in ${\bf R}^{n-1}$. Now for fixed $u$, then $x = s.u + r.v$ is in cylindrical polar coordinates with respect to $v$. We let $d\Omega(u,v)$ be the spherical volume form on the vectors perpendicular to $u$. Then as $u,v$ we have in terms of spherical polar coordinates
$$dV(x) = r^{n-3}d\Omega(u,v).dr.ds$$
 $$dV(x)dV(y) = (t-s)^{n-2}r^{n-3}d\Omega(u,v)d\Omega(u)dr.ds.d(t-s).$$
 Therefore
$$dV(x)dV(y)= (t-s)^{n-2}r^{n-3}d\Omega(u,v)d\Omega(u)dr.ds.dt.$$

We consider the vertical hyperbolic plane $P$ containing $x, y$. Then we see that $P$ intersects both  $P_0, P_a$ in semicircles centered about $w$. By definition of $r = |w|$ we have that the semicircles have radii
$r_1 = \sqrt{1-r^2}, r_2 = \sqrt{a^2-r^2}$ respectively. The geodesic with endpoints $x,y$ is contained in $P$. Therefore restricting to $P$ and scaling by $1/r_1$ we have that the length $L_{a}(x,y)$ is given by the previous lemma
$$L_a(x,y) = L_{\frac{r_2}{r_1}}\left(\frac{s}{r_1},\frac{t}{r_1}\right) = \frac{1}{2}\log\left(\frac{(t^2-r_1^2)(s^2-r_2^2)}{(t^2-r_2^2)(s^2-r_1^2)}\right).$$
 Then we have  $$F_n(l) = \frac{2}{V_{n-1}} \int_0^1 dr \int_{{\bf S}^{n-2}}d\Omega(u) \int_{{\bf S}^{n-3}}d\Omega(u,v)  \int_{|s| < r_1} ds\int_{t > r_2}  \frac{\log\left(\frac{(t^2-r_1^2)(s^2-r_2^2)}{(t^2-r_2^2)(s^2-r_1^2)}\right)}{(s-t)^{2n-2}}(t-s)^{n-2}r^{n-3}dt$$
Integrating over $u,v$ we get
$$F_n(l) = \frac{2V_{n-2}V_{n-3}}{V_{n-1}}  \int_0^1 r^{n-3}dr  \int^{r_1}_{- r_1} ds\int_{t > r_2}  \frac{\log\left(\frac{(t^2-r_1^2)(s^2-r_2^2)}{(t^2-r_2^2)(s^2-r_1^2)}\right)}{(t-s)^{n}}dt$$
We let $u =s/r_1, v = t/r_1$ and $b = r_2/r_1$ then 
$$F_n(l) = \frac{2V_{n-2}V_{n-3}}{V_{n-1}} \int_0^1 r^{n-3}dr \left(\frac{1}{r_1^{n-2}} \int^{1}_{-1} du\int^{\infty}_{b}  \frac{\log\left(\frac{(v^2-1)(u^2-b^2)}{(v^2-b^2)(u^2-1)}\right)}{(v-u)^{n}}dv\right)$$
We define
$$M_n(b) = \int^{1}_{-1} du\int^{\infty}_{b}  \frac{\log\left(\frac{(v^2-1)(u^2-b^2)}{(v^2-b^2)(u^2-1)}\right)}{(v-u)^{n}}dv$$
Therefore

\begin{equation}
F_n(l) =\frac{2V_{n-2}V_{n-3}}{V_{n-1}}\int_0^1 \frac{r^{n-3}}{\left(\sqrt{1-r^2}\right)^{n-2}}.M_n\left(\sqrt{\frac{a^2-r^2}{1-r^2}}\right)dr.
\label{Fn}
\end{equation}
We obtain an alternate integral  form  for $F_n$ by letting $x = \sqrt{\frac{a^2-r^2}{1-r^2}}$. then
$$x^2 = \frac{a^2-r^2}{1-r^2} \qquad r^2 = \frac{x^2-a^2}{x^2-1} \qquad 2rdr = \frac{(a^2-1)}{(x^2-1)^2}2xdx$$
Giving
$$r^{n-3}dr = r^{n-4}r.dr = \left(\frac{x^2-a^2}{x^2-1}\right)^{\frac{n}{2}-2}\left(\frac{a^2-1}{(x^2-1)^2}\right)xdx = 
 \frac{x(a^2-1)(x^2-a^2)^{\frac{n}{2}-2}}{(x^2-1)^{\frac{n}{2}}}dx$$ 
$$1-r^2= \frac{a^2-1}{x^2-1}\qquad \frac{1}{(\sqrt{1-r^2})^{n-2}} = \left(\frac{x^2-1}{a^2-1}\right)^{\frac{n}{2}-1}$$

Therefore
\begin{equation}
F_n(l) = \frac{2V_{n-2}V_{n-3}}{V_{n-1}} \int_a^\infty \left(\frac{x^2-a^2}{a^2-1}\right)^{\frac{n}{2}-2}\left(\frac{xM_n(x)}{x^2-1}\right)dx
\label{altFn}
\end{equation}
\subsection{Formula For $M_n$}
We give an explicit formula for $M_n$. in order to do so we first need to state some integral formulae. For  $n \geq 1$  we define the polynomial function $P_n$ by
$$P_n(z) = \sum_{k=1}^n \frac{x^k}{k}.$$
We also define $P_0(z) = 0$. We note that for $|x| < 1$, $P_n(x)$ is the first $n$ terms of the Taylor series of $-\log(1-x)$. 
We therefore define the function $L_n(x)$ by
$L_n(x) = \log|1-x| + P_n(x).$
For  $|x| < 1$ we have
$$L_n(x) = -\sum_{k=n+1}^\infty \frac{x^k}{k}.$$
We note that $L_0(x) = \log|1-x|$. We also note that $P_n(1) = 1 + \frac{1}{2}+\ldots+\frac{1}{n}$, the n$^{th}$ Harmonic number.

Then we have the following integral formula.
\begin{lemma} If $a \neq b \in{\bf R}$ and $n$ a positive integer then 
$$\int \frac{dx}{(x-a)(x-b)^n} \ =\  \frac{1}{(a-b)^n}\left( \log\left|\frac{x-a}{x-b}\right| + P_{n-1}\left(\frac{a-b}{x-b}\right)\right) = \frac{L_{n-1}\left(\frac{a-b}{x-b}\right)}{(a-b)^n}$$
\end{lemma}

{\bf Proof:} We note that
$$\frac{1}{x^n(x-c)} = \frac{1}{c^n}\left( \frac{1}{x-c} - \frac{x^n-c^n}{x^n(x-c)} \right)= \frac{1}{c^n}\left(\frac{1}{x-c} - \frac{Q_{n-1}(x)}{x^n}\right)$$
where
$$Q_{n-1}(x) = \frac{x^n-c^n}{x-c} = \sum_{k=1}^{n} x^{k-1}c^{n-k}$$
Therefore
$$\int \frac{1}{x^n(x-c)} dx = \frac{1}{c^n}\left(\int \frac{dx}{x-c} -  \int \frac{Q_{n-1}(x)}{x^n} dx\right) = \frac{1}{c^n}\left(\log|x-c| -  \int \frac{Q_{n-1}(x)}{x^n} dx\right) $$
As
$$ \frac{Q_{n-1}(x)}{x^n} =  \sum_{k=1}^{n} \frac{c^{n-k}}{x^{n-k+1}}$$
then
$$ \int \frac{Q_{n-1}(x)}{x^n}dx = \log|x| - \sum_{k=1}^{n-1} \frac{c^{n-k}}{(n-k).x^{n-k}} = \log|x|-P_{n-1}\left(\frac{c}{x}\right)$$
Therefore
$$\int \frac{1}{x^n(x-c)} dx  = \frac{1}{c^n}\left(\log\left|\frac{x-c}{x}\right|+P_{n-1}\left(\frac{c}{x}\right)\right) = \frac{L_{n-1}\left(\frac{c}{x}\right)}{c^n}$$
The result now follows by elementary substitution.
\eproof

Integrating by parts we get,
\begin{corollary} For $n \geq 2$
$$\int \frac{\log|x-a|}{(x-b)^n} dx = \frac{1}{n-1}\left(\frac{L_{n-2}\left(\frac{a-b}{x-b}\right)}{(a-b)^{n-1}} -\frac{\log|x-a|}{(x-b)^{n-1}}\right)$$
Furthermore for $k \geq 1$
$$\lim_{x \rightarrow a} \left(\frac{\log|x-a|}{(b-x)^{k}}-\frac{L_{n}\left(\frac{b-a}{b-x}\right)}{(b-a)^{k}}\right) = 
\frac{\log|b-a|-P_{n}(1)}{(b-a)^{k}}$$
\end{corollary}

\begin{lemma}
The function  $M_n:(1,\infty) \rightarrow {\bf R}_+$ has the explicit form
$$M_n(b) = \frac{1}{(n-1)(n-2)} \left(\frac{1}{(b-1)^{n-2}}\left( \log(\frac{(b+1)^2}{(4b)})+2P_{n-2}(1)-L_{n-3}(\frac{b-1}{b+1})-(-1)^n L_{n-3}(\frac{-b+1}{b+1})\right)+\right.$$
$$ \frac{1}{(b+1)^{n-2}}\left(-\log(\frac{(b-1)^2}{(4b)})-2P_{n-2}(1)+L_{n-3}(\frac{b+1}{b-1}) +(-1)^n L_{n-3}(\frac{-b-1}{b-1})\right)+$$
$$\left.\frac{1}{(2b)^{n-2}}\left(L_{n-3}(\frac{2b}{b+1})-L_{n-3}(\frac{2b}{b-1})\right)+\frac{1}{2^{n-2}}\left(L_{n-3}(\frac{2}{b+1})-(-1)^nL_{n-3}(\frac{-2}{b-1}\right)\right)$$
Furthermore $M_n$ satisfies 
$$\lim_{b \rightarrow 1^+} (b-1)^{n-2}M_n(b) = D_n = \frac{2P_{n-2}(1)}{(n-1)(n-2)} \qquad \mbox{and} \qquad \lim_{b \rightarrow \infty} \frac{b^{n-1}}{\log{b}}.M_n(b) = \frac{4}{n-1}.$$
\label{Mn}
\end{lemma}

{\bf Proof:}
We let $$I(u,b) = \int^{\infty}_{b}  \frac{\log\left(\frac{(v^2-1)(b^2-u^2)}{(v^2-b^2)(1-u^2)}\right)}{(v-u)^{n}}dv$$
then  $M_n(b) = \int_0^1 I(u,b) du$.
We split $I = I_1 + I_2$ where
$$I_1(u,b) =  \int^{\infty}_{b}  \frac{\log\left(\frac{b-u}{v-b}\right)}{(v-u)^{n}}dv\qquad I_2(u,b) = \int^{\infty}_{b}  \frac{\log\left(\frac{(v^2-1)(b+u)}{(v+b)(1-u^2)}\right)}{(v-u)^{n}}dv.$$
By the above lemma,
$$I_1(u,b) =\left. 
\frac{1}{n-1}\left(\frac{\log(\frac{v-b}{b-u})}{(v-u)^{n-1}}-\frac{L_{n-2}\left(\frac{b-u}{v-u}\right)}{(b-u)^{n-1}}\right)\right|_b^\infty$$
Therefore  by the above corollary
$$I_1(u,b)  = \frac{1}{n-1}\left(\frac{P_{n-2}(1)}{(b-u)^{n-1}}\right)$$

$$I_2(u,b) = \left. \left(\frac{-1}{n-1}\right)\frac{\log\left(\frac{(v^2-1)(b+u)}{(v+b)(1-u^2)}\right)}{(v-u)^{n-1}}\right|_b^\infty + \frac{1}{n-1}\int_b^\infty \frac{1}{(v-u)^{n-1}}\left(\frac{1}{v-1}+ \frac{1}{v+1}-\frac{1}{v+b}\right)dv$$
Evaluating we get
$$I_2(u,b) = \frac{1}{n-1}\left(\frac{\log\left(\frac{(b^2-1)(b+u)}{(2b)(1-u^2)}\right)}{(b-u)^{n-1}} + \int_b^\infty \frac{1}{(v-u)^{n-1}}\left(\frac{1}{v-1}+ \frac{1}{v+1}-\frac{1}{v+b}\right)dv\right)$$
Therefore
$$M_n(b) = \frac{1}{n-1}\left(\int_{-1}^{1} \frac{\log\left(\frac{(b^2-1)(b+u)}{(2b)(1-u^2)}\right)+P_{n-2}(1)}{(b-u)^{n-1}}du  +\int_{-1}^{1} \int_b^\infty \frac{1}{(v-u)^{n-1}}\left(\frac{1}{v-1}+ \frac{1}{v+1}-\frac{1}{v+b}\right)dv du\right)$$
We once again split this integral as follows;
$$J_1(b) = \frac{1}{n-1}\left( \int_{-1}^{1} \frac{\log\left(\frac{(b^2-1)}{(2b)}\right)+P_{n-2}(1)}{(b-u)^{n-1}}du\right)$$
$$J_2(b) =  \frac{1}{n-1}\left(\int_{-1}^{1} \frac{\log\left(\frac{b+u}{(1-u)(1+u)}\right)}{(b-u)^{n-1}}du\right)$$
$$J_3(b) =  \frac{1}{n-1}\left(\int_{-1}^{1} \int_b^\infty \frac{1}{(v-u)^{n-1}}\left(\frac{1}{v-1}+ \frac{1}{v+1}-\frac{1}{v+b}\right)dv du \right)$$
We now evaluate each of $J_1, J_2, J_3$.

{\bf $J_1$:} Integrating we get
$$J_1(b) =   \frac{\log\left(\frac{(b^2-1)}{(2b)}\right)+P_{n-2}(1)}{(n-1)(n-2)}\left(\frac{1}{(b-1)^{n-2}} - \frac{1}{(b+1)^{n-2}}\right)$$
{\bf $J_2$:} By the above lemma
$$J_2(b) = \left.\frac{1}{(n-1)(n-2)}\left(\frac{\log\left(\frac{b+u}{(1-u)(1+u)}\right)}{(b-u)^{n-2}} 
-\frac{L_{n-3}(\frac{2b}{b-u})}{(2b)^{n-2}}+\frac{L_{n-3}(\frac{b-1}{b-u})}{(b-1)^{n-2}}+\frac{L_{n-3}(\frac{b+1}{b-u})}{(b+1)^{n-2}} \right)\right|_{-1}^{1}$$
Therefore
$$(n-1)(n-2)J_2(b) = \left(\frac{\log(\frac{b+1}{2})-L_{n-3}(\frac{b-1}{b+1})}{(b-1)^{n-2}}\right)+ \left(\frac{L_{n-3}(\frac{2b}{b+1})-L_{n-3}(\frac{2b}{b-1}) }{(2b)^{n-2}}\right)+ \left(\frac{-\log(\frac{b-1}{2})+L_{n-3}(\frac{b+1}{b-1})}{(b+1)^{n-2}}\right)+ $$
$$ \lim_{u\rightarrow 1-} \left( -\frac{\log(1-u)}{(b-u)^{n-2}}+ \frac{L_{n-3}(\frac{b-1}{b-u})}{(b-1)^{n-2}}\right)+\lim_{u\rightarrow -1^+} \left( \frac{\log(1+u)}{(b-u)^{n-2}}- \frac{L_{n-3}(\frac{b+1}{b-u})}{(b+1)^{n-2}}\right)$$
By the corollary above
$$ \lim_{u\rightarrow 1-} \left( -\frac{\log(1-u)}{(b-u)^{n-2}}+ \frac{L_{n-3}(\frac{b-1}{b-u})}{(b-1)^{n-2}}\right) = 
\frac{P_{n-3}(1) -\log(b-1)}{(b-1)^{n-2}}$$
$$\lim_{u\rightarrow -1^+} \left( \frac{\log(1+u)}{(b-u)^{n-2}}- \frac{L_{n-3}(\frac{b+1}{b-u})}{(b+1)^{n-1}}\right) = 
\frac{\log(b+1)-P_{n-3}(1) }{(b+1)^{n-2}}$$
Therefore
$$(n-1)(n-2)J_2(b) = \left(\frac{\log(\frac{b+1}{2(b-1)})+ P_{n-3}(1)-L_{n-3}(\frac{b-1}{b+1})}{(b-1)^{n-2}}\right)+ \left(\frac{L_{n-3}(\frac{2b}{b+1})-L_{n-3}(\frac{2b}{b-1}) }{(2b)^{n-2}}\right)$$
$$+ \left(\frac{-\log(\frac{b-1}{2(b+1)}) -P_{n-3}(1) +L_{n-3}(\frac{b+1}{b-1})}{(b+1)^{n-2}}\right)$$

{\bf $J_3$:} Switching the order of integration we get
$$J_3(b) =   \frac{1}{(n-1)(n-2)}\int_{b}^{\infty} \left(\frac{1}{v-1}+ \frac{1}{v+1}-\frac{1}{v+b}\right)\left(\frac{1}{(v-1)^{n-2}}-\frac{1}{(v+1)^{n-2}}\right)dv$$
Therefore
$$J_3(b) = \frac{1}{(n-1)(n-2)}\left( \frac{\frac{1}{n-2}}{(b-1)^{n-2}} -\frac{\frac{1}{n-2}}{(b+1)^{n-2}}+\frac{L_{n-3}(\frac{2}{b+1})}{2^{n-2}}-\frac{L_{n-3}(\frac{-2}{b-1})}{(-2)^{n-2}}+\frac{L_{n-3}(\frac{-b-1}{b-1})}{(-b-1)^{n-2}}-\frac{L_{n-3}(\frac{-b+1}{b+1})}{(-b+1)^{n-2}}\right)$$

Combining $J_1, J_2, J_3$ we have
 $$M_n(b) = \frac{1}{(n-1)(n-2)} \left(\left( \frac{\log\left(\frac{(b+1)^2}{4b}\right)+P_{n-2}(1)+P_{n-3}(1)-L_{n-3}(\frac{b-1}{b+1}) +\frac{1}{n-2}-(-1)^n L_{n-3}(\frac{-b+1}{b+1})}{(b-1)^{n-2}}\right)+\right.$$
$$ \left(\frac{-\log\left(\frac{(b-1)^2}{(4b)}\right)-P_{n-2}(1)-P_{n-3}(1)+L_{n-3}(\frac{b+1}{b-1}) -\frac{1}{n-2}+(-1)^n L_{n-3}(\frac{-b-1}{b-1})}{(b+1)^{n-2}}\right)+$$
$$\left.\left(\frac{\frac{1}{b^{n-2}}(L_{n-3}(\frac{2b}{b+1})-L_{n-3}(\frac{2b}{b-1}))+L_{n-3}(\frac{2}{b+1})-(-1)^nL_{n-3}(\frac{-2}{b-1})}{2^{n-2}}\right)\right)$$
Noting that $P_{n-3}(1) + \frac{1}{n-2} = P_{n-2}(1)$ we simplify to get
$$M_n(b) = \frac{1}{(n-1)(n-2)} \left(\frac{1}{(b-1)^{n-2}}\left( \log(\frac{(b+1)^2}{(4b)})+2P_{n-2}(1)-L_{n-3}(\frac{b-1}{b+1})-(-1)^n L_{n-3}(\frac{-b+1}{b+1})\right)\right.$$
$$ +\frac{1}{(b+1)^{n-2}}\left(-\log(\frac{(b-1)^2}{(4b)})-2P_{n-2}(1)+L_{n-3}(\frac{b+1}{b-1}) +(-1)^n L_{n-3}(\frac{-b-1}{b-1})\right)$$
$$+\left.\frac{1}{(2b)^{n-2}}\left(L_{n-3}(\frac{2b}{b+1})-L_{n-3}(\frac{2b}{b-1})\right)+\frac{1}{2^{n-2}}\left(L_{n-3}(\frac{2}{b+1})-(-1)^nL_{n-3}(\frac{-2}{b-1}\right)\right)$$

It follows from the above directly that $M_n$ satisfies
$$\lim_{b\rightarrow 1^+} (b-1)^{n-2}M_n(b) = \frac{2P_{n-2}(1)}{(n-1)(n-2)}.$$
We now consider the limit as $b$ tends to infinity. As  $L_{n-3}$ is continuous except at $x = 1$, we replace the $L_{n-3}$ terms by the limiting value if the argument does not tend to 1 or is infinite. For other values we  substitute $P_n(x) = \log|1-x| + P_n(x)$ and collect terms. Therefore
$$\lim_{b \rightarrow \infty}(b-1)^{n-2}M_n(b) =  \frac{1}{(n-1)(n-2)} \lim_{b\rightarrow \infty} \left(\left( -\log(4)+2P_{n-2}(1)-\log(\frac{2}{b+1})- P_{n-3}(1)-(-1)^nL_{n-3}(-1)\right)\right.$$
$$ +\left(\log(4)-2P_{n-2}(1)+\log(\frac{2}{b-1})+P_{n-3}(1) +(-1)^nL_{n-3}(-1)\right)$$
$$+\left.\frac{1}{2^{n-2}}\left(L_{n-3}(2)-L_{n-3}(2)\right)+\frac{(b-1)^{n-2}}{2^{n-2}}\left(L_{n-3}(\frac{2}{b+1})-(-1)^nL_{n-3}(\frac{-2}{b-1}\right)\right)$$
After cancellation we have
$$\lim_{b\rightarrow 1^+} (b-1)^{n-2}M_n(b) = \frac{1}{(n-1)(n-2)} \lim_{b\rightarrow \infty} \frac{(b-1)^{n-2}}{2^{n-2}}\left(L_{n-3}(\frac{2}{b+1})-(-1)^nL_{n-3}(\frac{-2}{b-1}\right).$$
For $b$ large we can use the Taylor series expansion for $L_{n-3}$ to get  
$$L_{n-3}\left(\frac{2}{b+1}\right) =  \frac{1}{n-2}\left(\frac{2}{b+1}\right)^{n-2} + \epsilon_1.\left(\frac{2}{b+1}\right)^{n-3} \qquad L_{n-3}\left(\frac{-2}{b-1}\right) =  \frac{1}{n-2}\left(\frac{-2}{b-1}\right)^{n-2} + \epsilon_2.\left(\frac{-2}{b-1}\right)^{n-3}$$
Therefore 
$$ \lim_{b\rightarrow \infty} \frac{(b-1)^{n-2}}{2^{n-2}}\left(L_{n-3}(\frac{2}{b+1})-(-1)^nL_{n-3}(\frac{-2}{b-1})\right)= $$
$$\lim_{b\rightarrow \infty} \left(\frac{1}{n-2} + \epsilon_1\left(\frac{2}{b+1}\right)-(-1)^n\left(\frac{(-1)^{n-2}}{n-2} + \epsilon_2\frac{(-1)^{n-2}.-2}{b-1}\right) \right)= 0$$
We now use this fact to show the exact asyptotic behavior of $M_n(b)$ as $b$ tends to infinity. From the above we have   
$$\lim_{b \rightarrow \infty}b^{n-2}M_n(b) =0.$$
We now expand $M_n(b)$ around infinity by writing as a Taylor series in $1/b$. From the form of $M_n$ we must therefore have that for large $b$
$$M_n(b) = \frac{c_0+d_0\log{b}}{b^{n-1}} + \frac{c_1+ d_1\log{b}}{b^{n}} + \ldots$$
Thus we have that the leading terms is
$$M_n(b) \simeq \frac{d_0 \log{b}}{b^{n-1}} \qquad \mbox{ or } \qquad  \lim_{b \rightarrow \infty}b^{n-1}\log{b}.M_n(b) =d_0.$$
To find $d_0$ we need only gather the $\log{b}$ terms. We can replace terms of the form $\log(sb+t)$ by $\log{b}$ and also can ignore terms $L_{n-3}(m(b))$ where $m(b)$ does not tend to $1$ as $b$ tends to infinity. Therefore

$$M_n(b) \simeq \frac{1}{(n-1)(n-2)}\left(\frac{2\log{b}}{(b-1)^{n-2}}-\frac{2\log{b}}{(b+1)^{n-2}}\right) \simeq  \frac{1}{(n-1)(n-2)}\frac{2\log{b}}{b^{n-2}}\left(\frac{1}{(1-1/b)^{n-2}}-\frac{1}{(1+1/b)^{n-2}}\right)$$
Therefore
$$M_n(b) \simeq \frac{1}{(n-1)(n-2)}\frac{2\log{b}}{b^{n-2}}\left(1 + (n-2)\frac{1}{b}-\left(1-(n-2)\frac{1}{b}\right)\right) = \left(\frac{4}{n-1}\right)\frac{\log{b}}{b^{n-1}}.$$
\eproof

\begin{corollary}
For $n$ even there is a closed form for $F_n$.
\end{corollary}
{\bf Proof:} We note by equation \ref{altFn}, for all $n$ we have 
$$F_n(l) = \frac{2V_{n-2}V_{n-3}}{V_{n-1}} \int_a^\infty \left(\frac{x^2-a^2}{a^2-1}\right)^{\frac{n}{2}-2}\left(\frac{xM_n(x)}{x^2-1}\right)dx.$$
By the above lemma, $M_n$ is in terms of rational functions and rational functions times sums of linear log functions $\log(a+ bx)$. Therefore for $n$ even the function to be integrated in the above equation for $F_n$ is a sum  of rational functions and rational functions times linear log functions. As such functions have explicit formulae, we can therefore find formulae for $F_n$ for $n$ even.
As the number of terms of $M_n$ is linear in $n$ we obtain approximately a quadratic number of terms for $F_n$.
\eproof

\subsection{Three-dimensional case}
For $n =3$ we have
$$F_3(l) = 2.\int_0^1 \frac{M_3\left(\sqrt{\frac{a^2-r^2}{1-r^2}}\right)}{\sqrt{1-r^2}}.dr.$$
Using the above formula we have that 
$$M_3(x) = \frac{2}{x^2-1}(1-\log(2)) - \frac{1}{2x}\left(\frac{x-1}{x+1}\right)\log(x-1)+ \frac{1}{2x}\left(\frac{x+1}{x-1}\right)\log(x+1)$$
Also we have
$$M_3(x) \simeq \frac{1}{x-1}$$
for $x$ close to 1.

\subsection{Four-dimensional case}
For $n =4$ we have
$$F_4(l) = 8.\int_0^1 \frac{r.M_4\left(\sqrt{\frac{a^2-r^2}{1-r^2}}\right)}{1-r^2}.dr.$$
Using the above formula we have that 
$$M_4(b) = \frac{1}{6} \left(\frac{3 + 2\log(\frac{(b+1)^2}{4b})}{(b-1)^2}-\frac{3 + 2\log(\frac{(b-1)^2}{4b})}{(b+1)^2} +\frac{\log(\frac{b-1}{b+1})+ \frac{b}{b+1}-\frac{b}{b-1}}{2b^2} + \frac{\log(\frac{b-1}{b+1})+ \frac{1}{b+1} +\frac{1}{b-1}}{2}\right)
$$
Also we have
$$M_4(x) \simeq \frac{1}{2}\frac{1}{(x-1)^2}$$
for $x$ close to 1.

\section{Properties of $F_n$}

We  now describe the properties of $F_n$. In particular we complete the proof of Theorem 1.

\begin{lemma}
The function $F_n$  satisfies the following:
\begin{enumerate}
\item There exists a $C_n > 0$ such that
$$F_n(l) \leq \frac{C_n}{(e^l-1)^{n-2}}$$
\item $F_n$ is continuous monotonically decreasing 
\item 
$$\lim_{l \rightarrow 0}l^{n-2}F_n(l) = K_n = \frac{2\pi^{\frac{n-3}{2}}P_{n-2}(1)\Gamma(\frac{n}{2}+1)\Gamma(\frac{n}{2}-1)}{n.\Gamma(\frac{n+1}{2})\Gamma(n-1)}$$
\item 
$$\lim_{l\rightarrow \infty} \frac{e^{(n-1)l}}{l}. F_n(l) =  \frac{(n-2)\pi^{\frac{n-2}{2}}\Gamma(\frac{n}{2}-1)}{(\Gamma(\frac{n+1}{2}))^2}$$
\end{enumerate}
\label{Fnprop}
\end{lemma}

{\bf Proof:}
\begin{enumerate}
\item For $a = e^l$ we have
$$F_n(l) = \frac{2V_{n-2}V_{n-3}}{V_{n-1}}\int_0^1 \frac{r^{n-3}}{\left(\sqrt{1-r^2}\right)^{n-2}}.M_n\left(\sqrt{\frac{a^2-r^2}{1-r^2}}\right)dr$$
with
$$M_n(b) = \int_{-1}^{1} ds \int_{b}^{\infty}   \frac{\log\left(\frac{(t^2-1)(s^2-b^2)}{(t^2-b^2)(s^2-1)}\right)}{(t-s)^{n}}dt.$$
By lemma \ref{Mn} we have that there exists a $B_n>0$ such that
$$M_n(x) \leq \frac{B_n}{(x-1)^{n-2}}.$$
Therefore 
$$F_n(l) \leq \frac{2B_nV_{n-2}V_{n-3}}{V_{n-1}}\int_0^1 \frac{r^{n-3}}{\left(\sqrt{1-r^2}\right)^{n-2}.\left(\sqrt{\frac{a^2-r^2}{1-r^2}}-1\right)^{n-2}}dr.$$
$$F_n(l) \leq \frac{2B_nV_{n-2}V_{n-3}}{V_{n-1}}\int_0^1 \frac{r^{n-3}}{\left(\sqrt{a^2-r^2}-\sqrt{1-r^2}\right)^{n-2}}dr.$$
Rationalizing the denominator we get 
$$F_n(l) \leq \frac{2B_nV_{n-2}V_{n-3}}{V_{n-1}}\int_0^1 \frac{r^{n-3}(\sqrt{a^2-r^2}+\sqrt{1-r^2})^{n-2}}{(a^2-1)^{n-2}}dr.$$
Therefore $$(a-1)^{n-2}F_n(l) \leq \frac{2B_nV_{n-2}V_{n-3}}{V_{n-1}}\int_0^1 r^{n-3}\left(\frac{\sqrt{a^2-r^2}+\sqrt{1-r^2}}{a+1}\right)^{n-2}dr.$$
We have for $a \geq 0, 0 \leq r \leq 1$
$$0 \leq \frac{\sqrt{a^2-r^2}+\sqrt{1-r^2}}{a+1} \leq 1.$$
Therefore we have a $C_n > 0$ such that
$$(e^l-1)^{n-2}F_{n}(l) \leq C_n \qquad \mbox{ giving }\qquad F_n(l) \leq \frac{C_n}{(e^l-1)^{n-2}}.$$ 

\item We have for $a=e^l$ then by equation \ref{altFn},
$$F_n(l) = \frac{2V_{n-2}V_{n-3}}{V_{n-1}} \int_a^\infty \left(\frac{x^2-a^2}{a^2-1}\right)^{\frac{n}{2}-2}\left(\frac{xM_n(x)}{x^2-1}\right)dx.$$
The function $\frac{x^2-a^2}{a^2-1}$ is monotonically decreasing in $a$. Therefore for $n \geq 4$ then $\frac{n}{2}-2 \geq 0$ and if $a < b$
$$ \int_b^\infty \left(\frac{x^2-b^2}{b^2-1}\right)^{\frac{n}{2}-2}\left(\frac{xM_n(x)}{x^2-1}\right)dx \leq \int_b^\infty \left(\frac{x^2-a^2}{a^2-1}\right)^{\frac{n}{2}-2}\left(\frac{xM_n(x)}{x^2-1}\right)dx \leq  \int_a^\infty \left(\frac{x^2-a^2}{a^2-1}\right)^{\frac{n}{2}-2}\left(\frac{xM_n(x)}{x^2-1}\right)dx $$
Thus $F_n$ is monotonic for $n \geq 4$.
The case for $n=3$ we have an explicit form for $M_3$ given by
$$M_3(x) = \frac{2}{x^2-1}(1-\log(2)) - \frac{1}{2x}\left(\frac{x-1}{x+1}\right)\log(x-1)+ \frac{1}{2x}\left(\frac{x+1}{x-1}\right)\log(x+1).$$
This function is monotonically decreasing. As 
$$F_3(l) = 2.\int_0^1 \frac{M_3\left(\sqrt{\frac{a^2-r^2}{1-r^2}}\right)}{\sqrt{1-r^2}}.dr$$
then it follows that $F_3$ is also monotonically decreasing.

\item We now analyse the behavior of $F_n(l)$ as $l$ tends to zero. By lemma \ref{Mn} we have that
$$\lim_{b \rightarrow 1^+} (b-1)^{n-2}M_n(b) = D_n$$.
Let $\epsilon >0$. then on $[0,1-\epsilon]$ we have
$$\lim_{a \rightarrow 1^+} \left(\sqrt{\frac{a^2-r^2}{1-r^2}}-1\right)^{n-2}M_n\left(\sqrt{\frac{a^2-r^2}{1-r^2}}\right) = D_n\qquad \mbox{uniformly}.$$
Simplifying we have
$$\left(\sqrt{\frac{a^2-r^2}{1-r^2}}-1\right)^{n-2} =  \frac{(a^2-1)^{n-2}}{(\sqrt{a^2-r^2}+\sqrt{1-r^2})(\sqrt{1-r^2})^{n-2}}$$
Therefore we have
$$\lim_{a \rightarrow 1^+} (a-1)^{n-2}M_n\left(\sqrt{\frac{a^2-r^2}{1-r^2}}\right) = D_n.(1-r^2)^{n-2}\qquad \mbox{uniformly on } r \in [1,1-\epsilon].$$
Therefore
$$\lim_{a \rightarrow 1^+} (a-1)^{n-2}\int_0^{1-\epsilon}  \frac{r^{n-3}}{\left(\sqrt{1-r^2}\right)^{n-2}}.M_n\left(\sqrt{\frac{a^2-r^2}{1-r^2}}\right)dr = D_n. \int_0^{1-\epsilon}  r^{n-3}(\sqrt{1-r^2})^{n-2} dr$$

As $M_n(x) \leq \frac{B_n}{(x-1)^{n-2}}$ we have as above that
$$\lim_{a \rightarrow 1^+} (a-1)^{n-2} \int_{1-\epsilon}^{1} \frac{r^{n-3}}{\left(\sqrt{1-r^2}\right)^{n-2}}.M_n\left(\sqrt{\frac{a^2-r^2}{1-r^2}}\right)dr \leq B_n \int_{1-\epsilon}^1 r^{n-3}dr \leq B_n.\epsilon.$$
Also for $n \geq 3$ 
$$\int_{1-\epsilon}^1  r^{n-3}(\sqrt{1-r^2})^{n-2} dr \leq \epsilon.$$
Therefore
$$\limsup_{a \rightarrow 1^+} \left|(a-1)^{n-2}\int_0^{1}  \frac{r^{n-3}}{\left(\sqrt{1-r^2}\right)^{n-2}}.M_n\left(\sqrt{\frac{a^2-r^2}{1-r^2}}\right)dr - D_n. \int_0^{1}  r^{n-3}(\sqrt{1-r^2})^{n-2} dr\right| \leq (B_n+D_n)\epsilon.$$
As $\epsilon$ is arbitrary we have
$$\lim_{l\rightarrow 0^+} l^{n-2}F_n(l) = \lim_{a \rightarrow a^+} (a-1)^{n-2}F_n(l) = \frac{2V_{n-2}V_{n-3}D_n}{V_{n-1}}. \int_0^1  r^{n-3}(\sqrt{1-r^2})^{n-2} dr.$$
Substituting $x = r^2$, and letting $B$ be the Beta function then 
$$\int_0^1  r^{n-3}(\sqrt{1-r^2})^{n-2} dr = \frac{1}{2} \int_0^1 x^{\frac{n}{2}-2}(1-x)^{\frac{n}{2}-1} dx = \frac{1}{2}B(\frac{n}{2}-1,\frac{n}{2}) = \frac{\Gamma(\frac{n}{2}-1). \Gamma(\frac{n}{2})}{2.\Gamma(n-1)}$$
 Combining we get
 $$\lim_{l\rightarrow 0^+} l^{n-2}F_n(l) =  \frac{D_n.V_{n-2}V_{n-3}\Gamma(\frac{n}{2}-1). \Gamma(\frac{n}{2})}{V_{n-1}.\Gamma(n-1)}$$
As the volume formula for the $n-$sphere is
$$V_n = \frac{(n+1)\pi^{\frac{n+1}{2}}}{\Gamma(\frac{n+3}{2})}$$
$$\lim_{l\rightarrow 0^+} l^{n-2}F_n(l) =  K_n = \frac{2\pi^{\frac{n-3}{2}}P_{n-2}(1)\Gamma(\frac{n}{2}+1)\Gamma(\frac{n}{2}-1)}{n.\Gamma(\frac{n+1}{2})\Gamma(n-1)}$$
We calculate $K_n$ for some small values. Starting with $n =3$ the first 10  values of $K_n$ are
$$\frac{\pi}{2}, 1, \frac{11\pi^2}{192}, \frac{5\pi}{54},\frac{137\pi^3}{30720}, \frac{7\pi^2}{1125}, \frac{121\pi^4}{458752},\frac{761\pi^3}{2315250}, \frac{7129\pi^5}{566231040},\frac{1342\pi^4}{93767625},\ldots$$

\item We now consider $F_n(l)$ for $l$ large. 
 Then  for large $a$ as 
$(a^2-r^2)/(1-r^2) > a^2-1$ for $0 < r < 1$, and we have that 
$$\lim_{a 
\rightarrow \infty} \frac{
\left(\sqrt{\frac{a^2-r^2}{1-r^2}}\right)^{n-1}}{\log\left(\sqrt{\frac{a^2-r^2}{1-r^2}}\right)}M_n\left(\sqrt{\frac{a^2-r^2}{1-r^2}}\right)  = \frac{4}{n-1} \qquad \mbox{uniformly on }[0,1].$$
Therefore on the interval $[0,1-\epsilon]$ we have
$$\lim_{a \rightarrow \infty} \frac{a^{n-1}}{\log{a}}M_n\left(\sqrt{\frac{a^2-r^2}{1-r^2}}\right)  =  \frac{4}{n-1}.(\sqrt{1-r^2})^{n-1}  \qquad \mbox{uniformly on }[0,1-\epsilon].$$
Therefore
$$\lim_{a \rightarrow \infty} \frac{a^{n-1}}{\log{a}} \int_0^{1-\epsilon} \frac{r^{n-3}}{(\sqrt{1-r^2})^{n-2}}.M_n\left(\sqrt{\frac{a^2-r^2}{1-r^2}}\right)dr = \frac{4}{n-1}\int_0^{1-\epsilon}r^{n-3}\sqrt{1-r^2}dr$$
As this holds for all $\epsilon$, the improper integral exists and we have
$$\lim_{a \rightarrow \infty} \frac{a^{n-1}}{\log{a}} \int_0^{1} \frac{r^{n-3}}{(\sqrt{1-r^2})^{n-2}}.M_n\left(\sqrt{\frac{a^2-r^2}{1-r^2}}\right)dr =  \frac{4}{n-1}\int_0^{1}r^{n-3}\sqrt{1-r^2}dr =  \frac{\sqrt{\pi}\Gamma(\frac{n}{2}-1)}{\Gamma(\frac{n+1}{2}).(n-1)}$$
Therefore
$$\lim_{l \rightarrow \infty} \frac{e^{(n-1)l}}{l} F_n(l) = \frac{2V_{n-2}V_{n-3}}{V_{n-1}}. \frac{\sqrt{\pi}\Gamma(\frac{n}{2}-1)}{\Gamma(\frac{n+1}{2}).(n-1)} = 
\frac{(n-2)\pi^{\frac{n-2}{2}}\Gamma(\frac{n}{2}-1)}{(\Gamma(\frac{n+1}{2}))^2}$$

\end{enumerate}
\eproof
The above lemma completes the proof of Theorem 1.

\section{Volume Bounds}
As an application of the above, we now consider lower bounds on volume for hyperbolic manifolds with totally geodesic boundary. Immediately from the monotonicity of $F_n$ we have the following;

\begin{lemma}
If $M$ is a finite volume hyperbolic manifold with totally  geodesic boundary with shortest orthogeodesic of length $l$ then
$$Vol_n(M) \geq F_n(l).$$
\end{lemma}

In the paper \cite{Miy94}, Miyamoto studies volumes of hyperbolic manifolds with totally geodesic boundary using the length shortest return path from the boundary to itself, which in terms of orthogeodesics,  is twice the length  shortest orthogeodesic. Using packing methods for neighborhoods of boundary components he obtains the following result.

\begin{theorem}{(Miyamoto \cite{Miy94})}
There exists a function $\rho_n:{\bf R}_+ \rightarrow {\bf R}_+$ such that if $M$ is a finite volume hyperbolic manifold with totally  geodesic boundary with shortest orthogeodesic of length $l$ then
$$Vol_n(M) \geq \rho_n(l).Vol_{n-1}(\partial M)$$
with equality  if and only if M is decomposed into truncated regular simplices of edge-length $2l$. 
Furthermore the function $\rho_n$ is monotonically increasing with $\rho_n(0) > 0$ giving
$$Vol_n(M) \geq \rho_n(0).Vol_{n-1}(\partial M).$$
\end{theorem}

We note that these bounds are quite different as functions of $l$ as $F_n(l)$ tends to infinity as $l\rightarrow 0$ while $\rho_n(l)$ tends to a finite value. Also the lower bound of Miyamoto is proportional to the volume of the boundary while the bound given by $F_n$ gives no obvious relation.

We now use the asymptotic behavior of $F_n$ to get bounds on the volume of a hyperbolic n-manifold with geodesic boundary in terms of the volume of the boundary similar to the bound of Miyamoto. We first describe the volume of an $r$-neighborhood of a boundary component in terms of $r$. 

\begin{lemma}
Let $X$ be a finite volume subset of ${\bf H}^{n-1} \subseteq {\bf H}^n$. Let $N^+_r(X)$ be  a one sided  $r$-neighborhood of $X$ in ${\bf H}^n$. Then
$$\mbox{Vol}_{n}(N^+_r(X)) = \mbox{Vol}_{n-1}(X).S_n(r)$$
where 
$$S_n(r) = \int_0^r \cosh^{n-1}(x)\ dx.$$
\label{nbd}
\end{lemma}
{\bf Proof:}
We consider $X$ on the vertical hyperplane $P$ in the upper half-space model of ${\bf H}^n$ given by $x_1 = 0$. Then by elementary hyperbolic geometry the boundary of $N_r^+(P)$ is the Euclidean plane $x_n = \tan(\theta_0).x_1$ for $\theta_0 $ satisfying  $\sinh{r} = \tan{\theta_0}$.

We let $x_1 = r\sin{\theta},  x_n = r\cos{\theta}$.  Then $(r,\theta,x_2,\ldots,x_{n-1}) \in N_r^+(X)$ if and only if
$(0,x_2,x_3,\ldots,r) \in X$ and $0 < \theta <  \theta_0$. Changing variables we get
$$\mbox{Vol}_{n}(N^+_r(X)) = \int_{N^+_r(X)} \frac{dx_1dx_2\ldots dx_n}{x^n_n} =  \int_{N^+_r(X)} \frac{rd\theta drdx_2\ldots dx_{n-1}}{r^n\cos^n{\theta}} =$$
$$ \int_X \frac{dx_2\ldots dx_{n-1}dr}{r^{n-1}}.\int_0^{\theta_0}
\frac{d\theta}{\cos^n{\theta}}.$$
As $x_n = r$ on $P$ we get 
$$\mbox{Vol}_{n}(N^+_r(X))  = \int_X \frac{dx_2\ldots dx_{n-1}dx_n}{x_n^{n-1}}.\int_0^{\theta_0}
\frac{d\theta}{\cos^n{\theta}} = \mbox{Vol}_{n-1}(X) .\int_0^{\theta_0}
\frac{d\theta}{\cos^n{\theta}}.$$
Therefore 
$$\mbox{Vol}_{n}(N^+_r(X))  =  \mbox{Vol}_{n-1}(X). S_n(r)$$
where
$$S_n(r) = \int_0^{\theta_0} \sec^n{\theta}\  d\theta$$
for $\sinh{r} = \tan{\theta_0}$. We let $\sinh{x} = \tan{\theta}$ then $\cosh{x} dx = \sec^2{\theta}d\theta$.
As $\cosh{x} = \sec{\theta}$ we have $dx = \sec{\theta}d\theta$ and
$$S_n(r) = \int_0^{\theta_0} \sec^n{\theta}\  d\theta = \int_0^{\theta_0} \sec^{n-1}(\theta) \sec{\theta}d\theta  = \int_0^r \cosh^{n-1}(x) dx.$$
\eproof

We now prove Theorem 2. We first restate the theorem.\newline
{\bf Theorem 2.} {\em  For $n \geq 3$, there exists a monotonically increasing function $H_n:{\bf R}_+ \rightarrow{\bf R}_+$ and a constant $C_n > 0$ such that if $M$ is a hyperbolic n-manifold with totally geodesic boundary  of area $A$ then 
$$\mbox{Vol}_n(M) \geq H_n(A) \geq C_n.A^{\left(\frac{n-2}{n-1}\right)}$$
}

{\bf Proof:}
We let $A = \mbox{Vol}_{n-1}(\partial M)$. We consider taking $\delta$ neighborhoods of $\partial M$ in $M$ denoted $N_{\delta}(\partial M)$. We let $\delta_0 > 0$ be the largest values such that the interior of $N_{\delta}(\partial M)$ is embedded. Then $M$ has an orthogeodesic $\alpha$ of length $d \leq 2\delta_0$ . Therefore there are two contributions to the volume, one part $V_1$ from the volume of the region associated with the orthogeodesic $\alpha$,  and the other $V_2$ from the neighborhood of the boundary with $\mbox{Vol}_{n}(M) \geq \max(V_1, V_2).$ Then we have
 $V_1 = F_n(d)$ and as $F_n$ is monotonically decreasing $V_1 \geq F_n(2\delta_0)$. As $V_2$ is the (one-sided) $\delta_0$ neighborhood of $\partial M$, by lemma \ref{nbd} we have 
$$V_2 = A.S_n(\delta_0).$$
We note that for small $x$, we have the simple approximation $S_n(x) \simeq x$, with $S_n(x) \geq \sinh{x}$ for all $x  > 0$. Also by the formula for $S_n$, we have that $S_n$ is convex.
  
\begin{figure}[h] 
   \centering
   \includegraphics[width=4in]{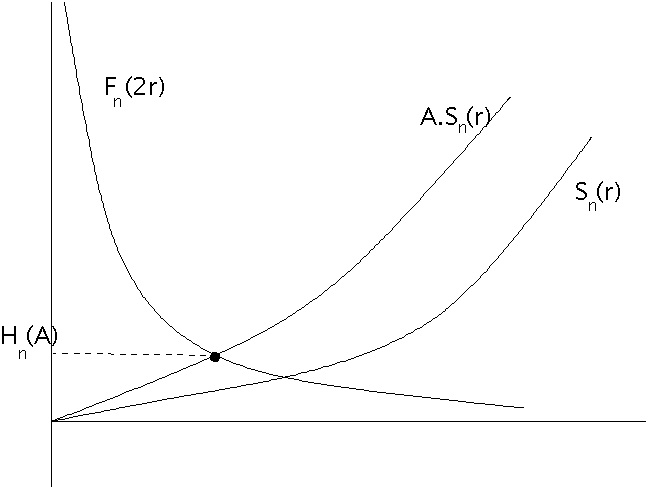} 
   \caption{Graph of $G_n(r)$}
\label{Gn}
\end{figure}

By lemma \ref{Fnprop}, we have that $F_n$ is a continuous monotonically decreasing positive function that tends to zero at infinity and tends to infinity at $0$. 
 We consider the continuous positive function $G_n(r) = \max(F_n(2r), A.S_n(r))$. Then $G_n$  has a positive minimum value depending only on $n$ and $A$ denoted $H_n(A)$. As $F_n$ is monotonically decreasing and $S_n$ monotonically increasing, we have that $H_n = J^{-1}$, where $J$ is the function
 $$J(x) = \frac{F_n(2x)}{S_n(x)}.$$
 We also note that $H_n(A)$ is the common value of the functions $F_n(2r), A.S_n(r)$ at their unique intersection point (see figure \ref{Gn}).
  
  Therefore we have $\mbox{Vol}_{n}(M) \geq G_n(\delta_0) \geq H_n(A)$. 
Obviously as $V_2$ is monotonic in volume, the function $H_n(A)$ is monotonic increasing in $A$ (see figure \ref{Gn}).

By lemma \ref{Fnprop}, there is a $K_n > 0$ such that
$$\lim_{l\rightarrow 0^+} l^{n-2}F_n(l) = K_n.$$ 
Therefore there is a $k_n$ such that for $l < 1$
$$F_n(r) \geq \frac{k_n}{r^{n-2}} = f_n(r).$$
 Also we have that $S_n(r) \geq A.r = s_n(r)$ for all $r$.
We consider the equation $f_n(2x) = s_n(x)$. This has a unique solution 
$$ \frac{k_n}{2^{n-2}.x^{n-2}} = A.x\qquad x = \left(\frac{k_n}{2^{n-2}A}\right)^{1/(n-1)} = \frac{k'_n}{A^{1/(n-1)}}$$ 
If $x < .5$ then we have
for $r \leq x$ then  $F_n(2r) \geq f_n(2r) \geq f_n(2x)$. Therefore $G_n(r) \geq f_n(2x)$ for $r \leq x$. 
Also for $r > x$ then $S_n(r) \geq s_n(r) \geq s_n(x)$. Therefore $G_n(r) \geq s_n(x)$ for $r > x$. 
Thus as $f_n(2x)  = s_n(x)$ then $G_n(r) \geq s_n(x)$. Thus $H_n(A) \geq s_n(x)$ giving
$$H_n(A) \geq A.x  = A. \frac{k'_n}{A^{1/(n-1)}} =  k'_n A^{\frac{n-2}{n-1}}$$
If $x \geq .5$ then  we have by monotonicity of $f_n, s_n$, that $f_n(1) \geq s_n(.5)$.
Therefore for $r \geq .5$ $S_n(r) \geq s_n(r) \geq s_n(.5)$. Therefore $G_n(r) \geq s_n(.5)$ for $r \geq .5$.
 If $r < .5$  then $F_n(2r) \geq f_n(2r) \geq f_n(1) \geq s_n(.5)$. Therefore $G_n(r) \geq s_n(.5)$ for $r \leq .5$. Thus $G_n(r) \geq s_n(.5)$ for all $r$. Thus $H_n(A) \geq s_n(.5) = A/2$ giving
 $$H_n(A) \geq \frac{1}{2}.A = C_n.A$$
 \eproof

We note that using the approximation $F_n(x) = \frac{K_n}{x^{n-2}}$ for $x$ small we get the approximate bound of
$$H_n(A) \stackrel{\sim}{\geq} \left(\frac{K_n.A}{2}\right)^{\frac{n-2}{n-1}}.$$

  \begin{figure} 
   \centering
   \includegraphics[width=4in]{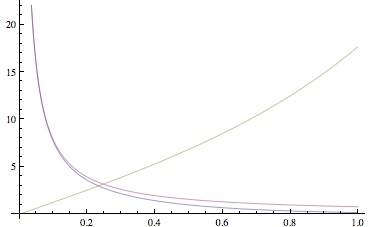} 
   \caption{Graph of $F_3(2x)$, $\frac{\pi}{4x}$ and $4\pi S_3(x)$}
\label{f3}
\end{figure}

\section{General lower bound Examples}
As the function $H_n(A)$ is monotonically increasing in $A$, if we have  a general lower bound on the volume of a closed hyperbolic $n-$manifold then we obtain a lower bound for a hyperbolic $(n+1)-$manifold with geodesic boundary.

{\bf Hyperbolic 3-manifolds}\newline
In \cite{KM91}, Kojima and Miyamoto showed that the lowest volume hyperbolic three manifold with totally geodesic boundary has boundary a genus two surface and volume 6.452.  Therefore $6.452$ is the best general lower bound. 

By the above lemma \ref{Fnprop}, $F_3(x) \simeq \frac{\pi}{2x}$ for small $x$. Graphing $F_3$ we see that this approximation is very close (see figure \ref{f3}).

For comparison with Kojima and Miyamoto, we note that as a closed hyperbolic surface has area at least $4\pi$, our methods give a lower bound of $H_3(4\pi)$ for the volume of a hyperbolic three-manifold with totally geodesic boundary.  Plotting $F_3(2x)$ and $4\pi.S_3(x)$ we see that $H_3(4\pi) \approx 2.986$ (see figure \ref{f3}).

 \end{document}